\documentclass[notitlepage, longbibliography]{revtex4-1}%
\usepackage{placeins}
\usepackage{pbox}
\usepackage{color}
\usepackage{relsize}
\usepackage{graphicx}
\usepackage[intlimits]{amsmath}
\usepackage{amsxtra, amssymb,fancyhdr, amsthm,latexsym}
\usepackage{MathShorthand}
\usepackage{cases}
\usepackage{verbatim}
\usepackage{float}
\usepackage{textcomp}
\usepackage[margin=0.75in]{geometry}
\usepackage{amsfonts}
\usepackage{amssymb}%
\providecommand{\U}[1]{\protect\rule{.1in}{.1in}}
\newtheorem{theorem}{Theorem}[section]

\newtheorem{prop}[theorem]{Proposition}

\begin{document}
\begin{empty}
\title{Metastable states for an aggregation model with noise}
\author{Joep H.M.~Evers$^{\dagger}$ and Theodore Kolokolnikov$^{\dagger}$}
\affiliation{$^{\dagger}$Department of Mathematics and Statistics, Dalhousie University,
Halifax, Canada}
\begin{abstract}
We study the long-time effect of noise on pattern formation for the aggregation
model. We consider aggregation kernels that generate
patterns consisting of two delta-concentrations.
Without noise, there is a one-parameter
family of admissible equilibria that
consist of two concentrations whose mass is not necessary equal.
We show that when a small amount of noise is added, the heavier
concentration ``leaks''
its mass towards the lighter concentration over a very long time
scale, eventually resulting in the equilibration of the two masses. We
use exponentially small asymptotics to derive
the long-time ODE's that quantify this mass exchange.
Our theory is validated using full numerical simulations of the
original model -- both of the original
stochastic particle system and its PDE limit.
Our formal
computations show that adding noise destroys the degeneracy in
the equilibrium solution and leads to a unique symmetric steady state.
\end{abstract}
\maketitle
\end{empty}


\section{Introduction}

Aggregation is an ubiquitous natural phenomenon that pervades both the animal
world and many inanimate physical systems. In the animal kingdom, group
formation is observed across all levels from bacterial colonies and insect
swarms to complex predator-prey interactions in fish, birds and mammals.
Aggregation is also present in physical systems of all scales from the
smallest (Bose-Einstein Condensates, DNA buckyball molecules, fluid vortices)
to the largest (galaxies). The emergence of group behaviour is often a
consequence of individuals (or atoms) following very simple rules, without any
external coordination.

\begin{empty}
\begin{figure}[t]
\centering
\includegraphics[width=\textwidth]{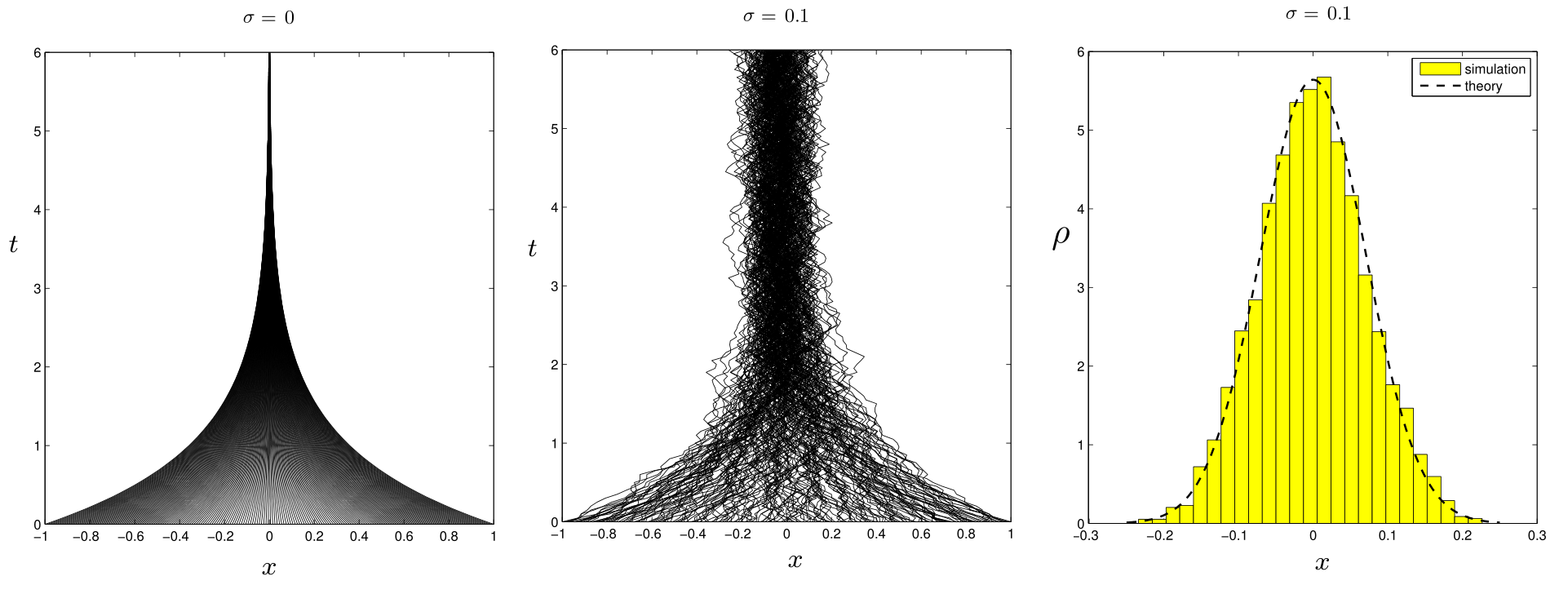}
\caption{LEFT: Particle simulation of deterministic system
(\ref{deterministic}) with $f(x)=-x$. Initial
conditions consist of $n=200$ particles equally spaced from each other.
Particles collapse into a single delta-concentration. MIDDLE: Same ODE system
but with noise added (\ref{noise}). The delta-concentration is ``diffused'' into
a Gaussian. RIGHT: Particle density distribution from the simulation in
the middle panel, averaged over 100
time-steps with $t=10.$ The dashed line is the predicted Gaussian
profile given by (\ref{rhogauss}).
}
\label{fig:simple}%
\end{figure}
\end{empty}

One of the simplest models that achieves aggregation is the so-called
aggregation model, which has been the subject of intense study in the last two
decades; refer to survey papers \cite{M&K, TBL, morale2005interacting,
bernoff2013nonlocal} and references therein. Mathematically, this model may be
written as a system of ODE's for $n$ particles%

\begin{equation}
\frac{dx_{j}}{dt}=\frac{1}{n}\sum_{k=1}^{n}f(x_{j}-x_{k})
\label{deterministic}%
\end{equation}
where the pairwise interaction force $f(x)$ is assumed to be the gradient of a
radial potential function, $f(x)=-\nabla P(\left\vert x\right\vert
)=-P^{\prime}(\left\vert x\right\vert )\frac{x}{\left\vert x\right\vert }$.
The strength of the force $f(x_{j}-x_{k})$ depends only on the distance
between the two particles $x_{j}$ and $x_{k}$, and it acts in the direction
between these particles. The system (\ref{deterministic})\ corresponds to
applying the method of steepest descent to determine the minimizer of the
pairwise-interaction energy,%
\[
E=\sum_{k,j}P(\left\vert x_{j}-x_{k}\right\vert ).
\]

\begin{empty}
\begin{figure}[t]
\centering
\includegraphics[width=\textwidth]{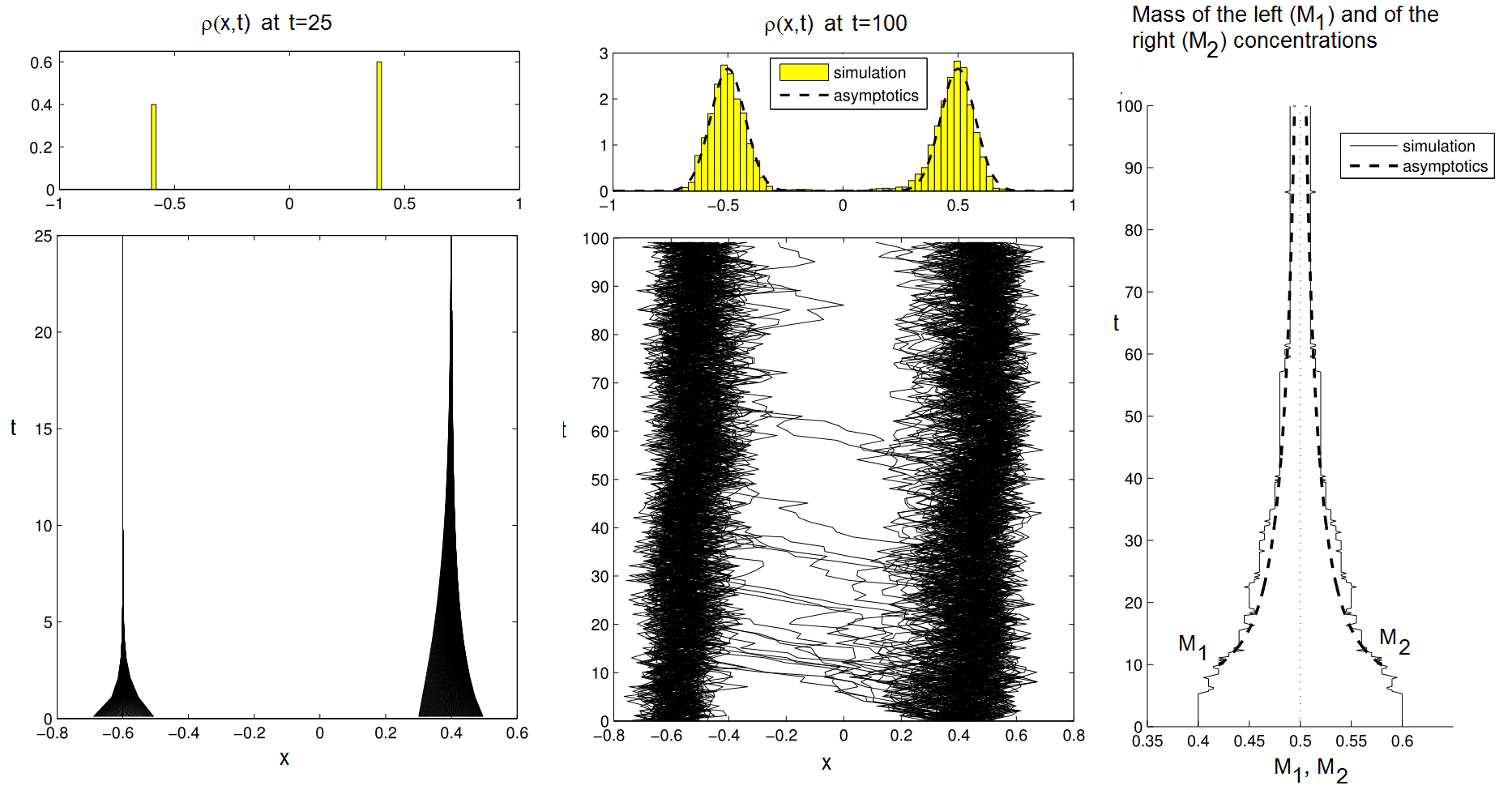}
\caption{
LEFT: Simulation of the deterministic system (\ref{deterministic}) with $f(x)=x-x^3$.
Initial conditions consist of 80 particles near $x=-0.5$ and 120 particles
near $x=0.5$ (corresponding to $M_1=0.4, M_2=0.6$). The long-time dynamics
approach two unequal delta-concentrations. The insert in top-left
shows the histogram of the final steady state.
MIDDLE: Same as the left figure, but with $\sigma=0.075$ noise added. The
concentration at the right initially has a larger mass, and it very
slowly leaks its mass towards a lighter concentration on the left.
RIGHT: The total mass of left and right concentrations of the
simulation in the middle is plotted. After a short transient period, a
slow mass exchange is apparent, with the two masses gradually
equilibrating. The dashed line denotes the asymptotic prediction given by
\eqref{odecubic}.
}
\label{fig:cubic}%
\end{figure}
\end{empty}

To get confinement, it is further assumed that the particles repel each other
at short distances and attract each other at longer distances (see also
\cite{Yao, Canizo, Simione}). In many cases this leads to the formation of
swarms. The assumption of long-range attraction and short-range repulsion
corresponds to $P(r)$ having a minimum at $r=r_{0}$ so that particles at a
distance bigger than $r_{0}$ are attracted to each other and those at distance
less than $r_{0}$ are repelling. These simple assumptions can give rise to
surprisingly complex steady states \cite{kolokolnikov2011stability,
von2012predicting, von2012soccer, fetecau2011swarm, bertozzi2012aggregation,
kolokolnikov2013singular, balague2013dimensionality, fellner2010stable,
von2013swarming} including \textquotedblleft soccer balls\textquotedblright%
\ in two and three dimensions \cite{kolokolnikov2011stability,
von2012predicting, kolokolnikov2013singular, balague2013dimensionality} as
well as steady states concentrating on points, curves and surfaces
\cite{fellner2010stable, kolokolnikov2013singular, balague2013dimensionality,
von2013swarming, hackett2012aggregation}. Of particular importance for the
dimensionality of the steady state is the strength of repulsion near the
origin \cite{kolokolnikov2011stability, von2012predicting,
balague2013dimensionality}. For this paper, we focus on the simplest case,
where the steady state concentrates on a finite number of points
(delta-concentrations), which can occur when the repulsion is sufficiently
weak at the origin.

In this work we are interested in how the noise that is inherently present in
most of the physical systems can affect the resulting steady state. That is,
we consider the model (\ref{deterministic})\ with noise, so that
(\ref{deterministic})\ is replaced by stochastic ODE's
\begin{equation}
dx_{j}=\frac{1}{n}\sum_{k=1}^{n}f(x_{j}-x_{k})dt+\sigma\sqrt{dt}%
\mathcal{N}_{j}; \label{noise}%
\end{equation}
here $\sigma\sqrt{dt}\mathcal{N}_{j}$ is the standard Wiener process 
with
standard deviation $\sigma$
; $\mathcal{N}_{j}$ denotes the standard normal
distribution of mean zero and variance 1. In the continuum limit as the number
of particles $n\rightarrow\infty,$ the average particle density distribution
$\rho$ is well approximated by the PDE \cite{morale2005interacting},
\begin{equation}
\rho_{t}+\nabla\cdot\left(  v\rho\right)  =\varepsilon^{2}\Delta
\rho,\ \ \ \ v=\int f\left(  x-y\right)  \rho\left(  y\right)  dy. \label{pde}%
\end{equation}
where $\varepsilon^{2}=\sigma^{2}/2.$ Equation (\ref{pde})\ is the starting
point for this paper.

The presence of noise can have a profound effect on the steady state,
especially if the steady state consists of point concentrations which can be
the case when the repulsion is sufficiently weak at the origin
\cite{hughes2013continuum}. As a motivating example, consider the simplest
case,
\begin{equation}
f(x)=-x.
\end{equation}
This corresponds to a pontential $P(r)=r^{2}/2$ that is purely attractive, and
is only weakly attractive at the origin (i.e. $P^{\prime}(0)=0$). Figure
\ref{fig:simple} shows the resulting one-dimensional simulations with and
without noise for the discrete system (\ref{noise}), as well as the associated
average density in the presence of noise -- computed by averaging the steady
state for the last 100 steps of the numerical simulation of (\ref{noise}) (for
simplicity, we used forward Euler method for these simulations). Without
noise, the particle density collapses to a single point\ (delta function). On
the other hand, the noise \textquotedblleft diffuses\textquotedblright\ the
delta function and the resulting average steady state density is a Gaussian,
\begin{equation}
\rho(x)=\frac{1}{\sqrt{2\pi\varepsilon^{2}}}\exp\left(  -\frac{x^{2}%
}{2\varepsilon^{2}}\right)  . \label{rhogauss}%
\end{equation}
as was already observed in \cite{hughes2013continuum}. Indeed, the steady
state satisfies $\left(  v\rho\right)  _{x}=\varepsilon^{2}\rho_{xx}$ where
$v(x)=\int-\left(  x-y\right)  \rho(y)dy=-x$, where we assumed that $\rho$ has
total mass $M=1$ and its centre of mass is $0$. Integrating once we get the
ODE $\varepsilon^{2}\rho_{x}=-x\rho$ whose solution is given by
(\ref{rhogauss}).

In this paper, we are interested in the effect of small amount of diffusion
when the steady state consists of more than one delta function. This
occurs when repulsion is present (so that the particles do not all collapse
into a single point), but the repulsion is sufficiently weak so that far-field
attraction causes \textquotedblleft clumping\textquotedblright\ into two or
more delta-concentrations. As was shown in \cite{fellner2010stable} in one
dimension, (and extended in \cite{kolokolnikov2013singular} to two and three
dimensions), the necessary condition for this to happen is that $P^{\prime
}(0)=0,$ or equivalently, $f(x)\sim cx$ as $x\rightarrow0$ for some
\emph{positive }constant $c>0$. To illustrate this phenomenon as well as the
results of the paper, consider the simplest such case, namely the double-well
potential, $P(r)=-r^{2}/2+r^{4}/4,$ so that%
\begin{equation}
f(x)=x-x^{3}. \label{cubic}%
\end{equation}
Figure \ref{fig:cubic}(left) illustrates the behaviour of the deterministic
system (\ref{deterministic}) with $n=200$ particles, starting with initial
conditions consisting of 80 particles near $x=-0.5$ and 120 particles near
$x=0.5$. After some transient time, the system evolves into a steady state
consisting of two delta concentrations, with 40\% of the mass at the left
concentration and 60\% at the right. The distance between the two
concentrations is $x=1$ corresponding to the root of (\ref{cubic}), and this
is trivially seen to be a steady state of (\ref{deterministic}) since
$f(0)=f(1)=0$. Moreover, as was shown in \cite{fellner2010stable,
kolokolnikov2013singular}, such steady state is actually stable.

Now suppose there is a small amount of noise present in the system (say
$\sigma=0.075$) while keeping all other parameters and initial conditions as
in Figure \ref{fig:cubic}(left). The result is shown in Figure \ref{fig:cubic}%
(middle and right). Initially, the system quickly settles to a
two-concentration asymmetric steady state with roughly 40\% of mass on the
left and 60\% of mass on the right, except that the noise \textquotedblleft
diffuses\textquotedblright\ the delta concentrations, so that the particles
constantly jiggle around, and the average density is approximately a Gaussian.
However on a much longer time-scale, there is a very slow exchange of mass
that takes place between the left and right concentrations, so that the bigger
concentration slowly leaks its mass towards the smaller, until the two
concentrations eventually equilibrate. In other words, adding even a small
amount of noise eventually \textquotedblleft symmetrizes\textquotedblright%
\ the asymmetric steady state over a long time. Quantifying this very slow
exchange of mass is the goal of this paper.

There are similarities between our work and \cite{Geigant}. They consider the
combination of interactions and diffusion on a one-dimensional interval with
periodic boundary conditions. Their main focus is on the $O(1)$ stability of
steady states consisting of multiple peaks. By contrast, here we are
interested in the slow-time evolution (metastability) of these states. We will
illustrate our result now.

Our main finding describes the exchange of mass between the asymmetric
\textquotedblleft diffused\textquotedblright\ concentrations. The precise
statement, for a general kernel $f(x)$, is given in Proposition \ref{prop:ode}%
. To illustrate the result, consider the cubic kernel (\ref{cubic}), refer to
Figure \ref{fig:cubic}. Starting with arbitrary initial conditions, the system
converges, on an $O(1)$ timescale, to a solution consisting of two
concentrations. These concentrations will in general have unequal masses (that
depend on initial conditions), and their asymptotic profile is a Gaussian
spike whose variance depends on their relative masses $M_{1},M_{2}$. Once
these spikes form, there is a very slow equilibration process, whereby the
mass of a heavier spike leaks towards the lighter one. This process is
meta-stable, meaning that it takes an exponentially long time (in
$\varepsilon$) for the masses to equilibrate. Specializing Proposition
\ref{prop:ode} to (\ref{cubic}), this slow mass exchange is described
asymptotically by an ODE\bes\label{odecubic}%
\begin{equation}
\frac{d}{dt}M_{1}\sim F(M_{1},M_{2})-F(M_{2},M_{1});\ \ \ M_{1}+M_{2}=M
\label{odeM}%
\end{equation}
where $M$ is the total mass, $M_{1}$ and $M_{2}$ are the masses of the two
spikes, and%
\begin{equation}
F(M_{1},M_{2})=\frac{M_{2}}{2\pi}\sqrt{\frac{\left(  2M_{2}-M_{1}\right)  }%
{M}}\left(  2M_{1}-M_{2}\right)  \exp\left(  \frac{-M_{2}\left(  2M_{1}%
-M_{2}\right)  ^{3}+O(\varepsilon^{2})}{4M^{3}\varepsilon^{2}}\right)  .
\label{Fm1m2}%
\end{equation}
\ees Note that the $O(\varepsilon^{2})$ term inside the exponent in
(\ref{Fm1m2}) becomes an $O(1)$ premultiplier when the exponential is
expanded, so the constant in front of the exponential is not asymptotically
accurate. A higher-order expansion is necessary to obtain that correction, and
we do not attempt this here; it is an open problem to do so. However the
exponential decay dominates the dynamics and to leading order, the constant in
front of the exponential is of lower order when the solution is plotted on a
log-time scale. This is explained further in \S \ref{sec: ODE derivation} -- see
Figure \ref{fig:ddp} there.

\begin{figure}[ptb]
\centering
\includegraphics[width=\textwidth]{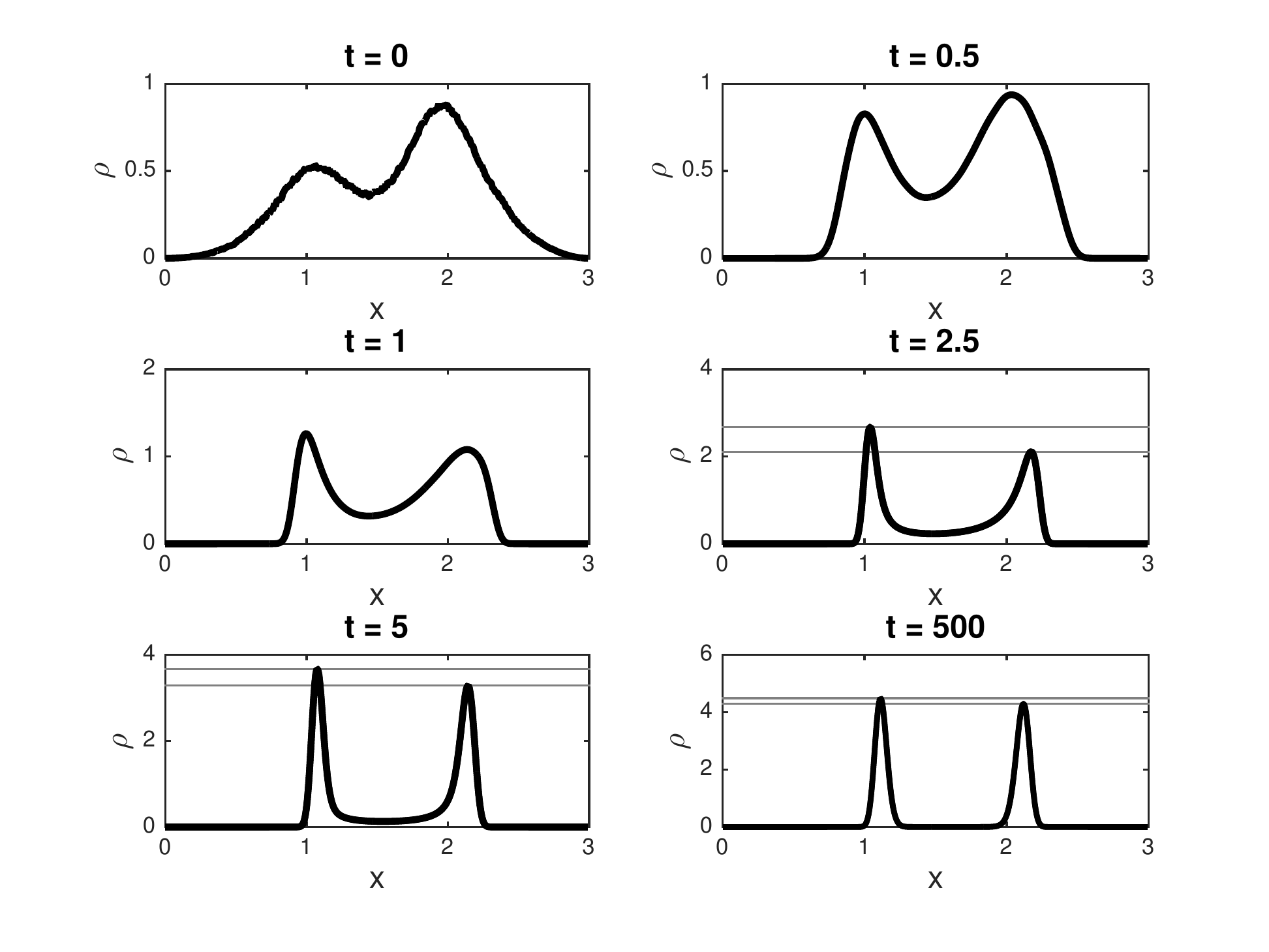}\caption{Evolution
of the density profile, for $f(x)=x(1-x^{2})$ and $\varepsilon^{2}=0.001$. In
the last three plots, the narrow grey lines indicate the heights of the two
maxima. Their heights converge at a slow timescale.}%
\label{fig:profile evol}%
\end{figure}

The summary of the paper is as follows. In \S \ref{sec:ss} we construct
asymptotically the quasi-steady state consisting of two Gaussians of masses
$M_{1},M_{2}.$ In \S \ref{sec:ode} we derive the equations of mass exchange
between $M_{1}$ and $M_{2}$, on an exponentially slow timescale, culminating
in Proposition \ref{prop:ode} which is the main result of this paper. From the
equations for mass exchange, we show that the masses equilibrate on a long
time scale, that is, $M_{1}=M_{2}$ is the unique global steady state of the
long-time dynamics. We compare the full numerical solution of the original
aggregation-diffusion PDE for several values of $\varepsilon$, to our asymptotic findings. We conclude with
some remarks in \S \ref{sec:discuss}.

\section{Quasi-steady state\label{sec:ss}}

In \cite{fellner2010stable, kolokolnikov2013singular} the authors constructed
a steady state of (\ref{pde}) with zero diffusion consisting of discrete
number $(N\geq2$)\ of delta-concentrations. This happens when $f(x)$ is linear
near the origin, such as for example (\ref{cubic}). More generally, we will
assume that:%
\begin{equation}
\left\{
\begin{tabular}
[c]{l}%
$f(x)$ is odd; \\
$f(x)$ has a positive root at $x=a;$ \\
$f(x)$ is $C^{1}$ at $x=0$ and $x=a$ with $f^{\prime}(0)>0$ and $f^{\prime}(a)<0.$%
\end{tabular}
\ \right.  \label{cond}%
\end{equation}
Under these assumptions, a two-delta steady state of (\ref{pde}) with
$\varepsilon=0$\ has the form%
\begin{equation}
\rho(x,t)\sim M_{1}\delta\left(  x-x_{1}\right)  +M_{2}\delta(x-x_{2}%
)\ \text{with }x_{2}-x_{1}=a,\ \ M_{1}+M_{2}=M \label{rhodelta}%
\end{equation}
where $M$ is the total mass. Upon substituting (\ref{rhodelta})\ into
$v=f\ast\rho$ we obtain that%
\begin{equation}
v(x)=M_{1}f(x-x_{1})+M_{2}f(x-x_{2}). \label{v}%
\end{equation}
Turning on $\varepsilon$ in (\ref{pde})\ \textquotedblleft
diffuses\textquotedblright\ the delta concentrations, so that the
$\delta\left(  x-x_{1}\right)  $ is replaced by a \textquotedblleft
spike\textquotedblright\ that has width of $O(\varepsilon)$. That is, we write%
\begin{equation}
\rho(x,t)\sim M_{1}\frac{1}{\varepsilon}w_{1}\left(  \frac{x-x_{1}%
}{\varepsilon}\right)  +M_{2}\frac{1}{\varepsilon}w_{2}\left(  \frac{x-x_{2}%
}{\varepsilon}\right)  \label{rho2}%
\end{equation}
where $w_{i}(y)$ is the spike profile that is to be computed, with
$\int_{-\infty}^{\infty}w_{i}(y)dy=1$ and $w_{i}(y)>0$ for all $y.$
Substituting (\ref{rho2})\ into $v=f\ast\rho$, and expanding the resulting
integral in terms of Taylor series, we find that the leading-order expression
for the velocity $v$\ is then still given by (\ref{v}) up to $O(\varepsilon)$
order. To compute the profile of the left spike $w_{1}(y),$ we let
$x=x_{1}+\varepsilon y$ and expand (\ref{rho2}). We have $v\left(
x_{1}\right)  =0$ so that%
\begin{equation}
v\left(  x_{1}+\varepsilon y\right)  \sim-\varepsilon c_{1}y,\ \ c_{1}%
=-v^{\prime}(x_{1})\sim-\left(  M_{1}f^{\prime}(0)+M_{2}f^{\prime}(a)\right)
. \label{vtaylor}%
\end{equation}
We then substitute (\ref{vtaylor})\ into the steady state equation after
discarding $\rho_{t}$ (this is the assumption that $\rho$ is a quasi-steady state).
Near $x=x_{1}+\varepsilon y$ we then obtain, up to exponentially small terms,%
\[
\left(  yc_{1}w_{1}\right)  _{y}\sim w_{1yy}.
\]
Assuming decay as $y\rightarrow\pm\infty$ yields%
\[
yc_{1}w_{1}+w_{1y}\sim0
\]
so that
\begin{equation}
w_{1}\sim\sqrt{\frac{c_{1}}{2\pi}}\exp\left(  \frac{-y^{2}}{2}c_{1}\right)  .
\end{equation}
A necessary condition for decay is that $c_{1}>0.$ Performing a similar
computation for $w_{2}$ we obtain the following result.

\begin{prop}
\label{prop:ss}Suppose that $f(x)$ satisfies conditions (\ref{cond}) and
suppose that $M_{1},M_{2}$ satisfy%
\begin{equation}
-f^{\prime}(0)/f^{\prime}(a)<M_{1}/M_{2}<-f^{\prime}(a)/f^{\prime}(0).
\label{mratio}%
\end{equation}
Then (\ref{pde})\ admits a quasi-equilibrium steady state that has the form
\begin{equation}
\rho(x,t)\sim\sum_{j=1}^{2}\frac{M_{j}}{\varepsilon}\sqrt{\frac{c_{j}}{2\pi}%
}\exp\left(  -\left(  \frac{x-x_{j}}{\varepsilon}\right)  ^{2}\frac{c_{j}}%
{2}\right)  \label{rhoinner}%
\end{equation}
with $x_{2}-x_{1}=a$ and $c_{j}=-v^{\prime}(x_{j})$ with $v$ given by
(\ref{v}); that is%
\begin{equation}
c_{1}=-\left(  M_{1}f^{\prime}(0)+M_{2}f^{\prime}(a)\right)  ;\ \ \ c_{2}%
=-\left(  M_{2}f^{\prime}(0)+M_{1}f^{\prime}(a)\right)  . \label{cj}%
\end{equation}
The masses $M_{1},M_{2}$ satisfy $M_{1}+M_{2}=M$ where $M$ is the total mass
that is determined by the initial conditions, $M=\int_{-\infty}^{\infty}%
\rho(x,0)dx.$
\end{prop}

\begin{figure}[ptb]
\centering
\includegraphics[width=0.5\textwidth,height=0.3\textwidth]{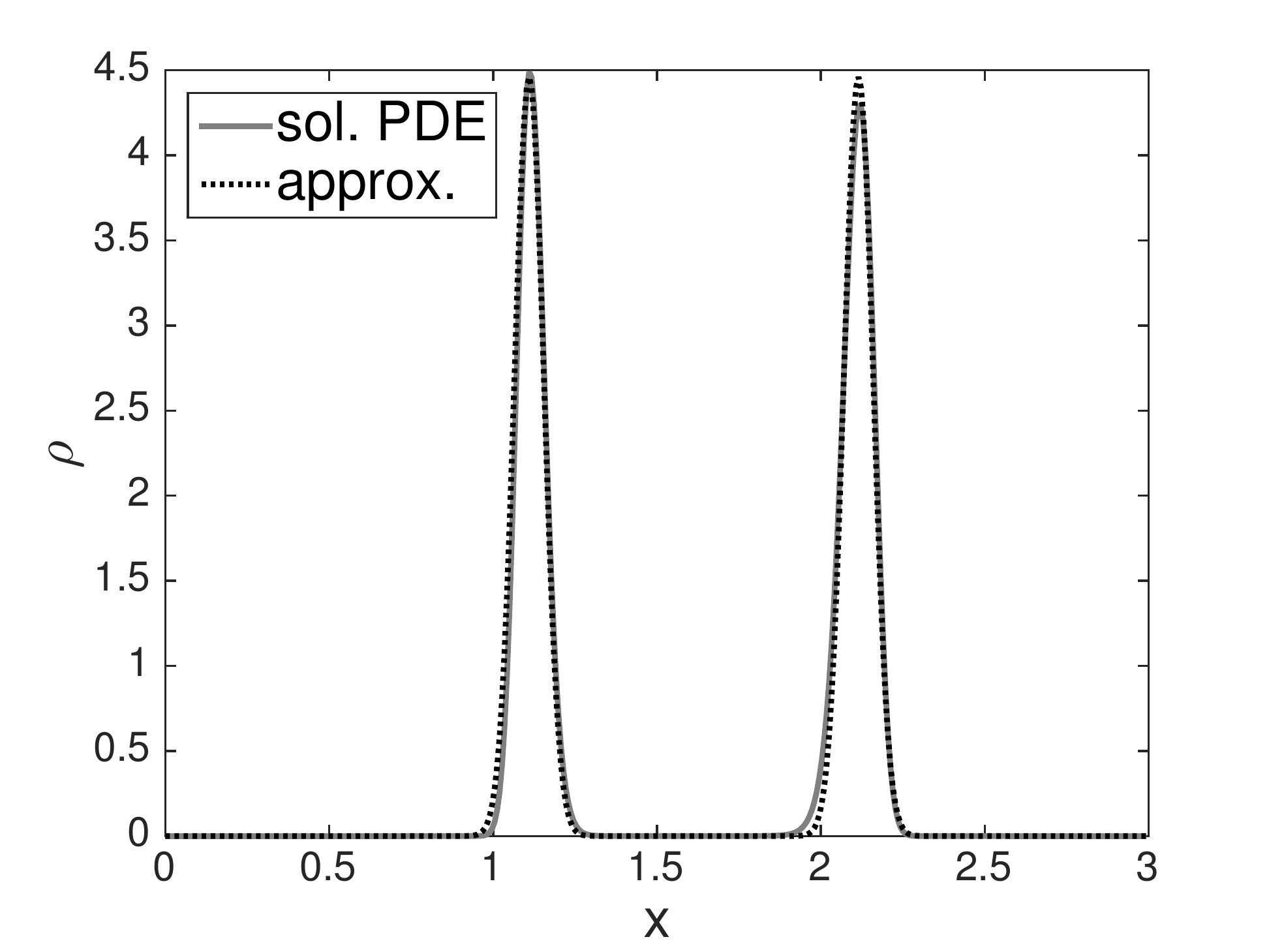}\caption{Density
profile at $t=500$, for $f(x)=x(1-x^{2})$ and $\varepsilon^{2}=0.001$; cf.
Figure \ref{fig:profile evol}. Superimposed is the approximation
\eqref{rhoinner} with $M_{1}=M_{2}$.}%
\label{fig:profile and approx}%
\end{figure}

As we show in \S \ref{sec:ode}, the masses $M_{1}$ and $M_{2}$ evolve on a
timescale much larger than the timescale at which this two-spike profile
forms. In Figure \ref{fig:profile and approx} we compare the long-time
solution of \eqref{pde}, as shown in the bottom right of Figure
\ref{fig:profile evol}, to the approximation \eqref{rhoinner} with
$M_{1}=M_{2}$.

Proposition \ref{prop:ss} generalizes naturally to $n$ concentrations. In this
case, the sum $\sum_{j=1}^{2}$ in (\ref{rhoinner})\ is replaced by $\sum
_{j=1}^{n}$, with%
\begin{equation}
c_{j}=-\sum_{k=1}^{n}M_{k}\,f^{\prime}(x_{j}-x_{k}).
\label{eqn: def eps j sq N}%
\end{equation}
and the condition $x_{2}-x_{1}=a$ is replaced by a system%
\begin{equation}
\sum_{k=1}^{n}M_{k}f(x_{j}-x_{k})=0,\ \ j=1\ldots n.
\end{equation}
Finally, the two conditions (\ref{mratio}) are replaced by $n$ conditions%
\begin{equation}
\sum_{k=1}^{n}M_{k}f^{\prime}(x_{j}-x_{k})<0.
\end{equation}

\section{Metastable dynamics\label{sec: ODE derivation}}

\label{sec:ode}We now derive the ODE describing the slow-time dynamics for the
mass exchange between $M_{1}$ and $M_{2}.$ The starting point is the PDE
(\ref{pde})\ with $v$ as given by (\ref{v}). Note from the expansion for $v$
near $x_{1}$ (\ref{vtaylor})\ that $v(x_{1})=0$ and $v^{\prime}(x_{1})<0$
(which implies that the concentration at $x=x_{1}$ attracts nearby points).
Similarly $v(x_{2})=0$ and $v^{\prime}(x_{2})<0.$ By continuity of $v(x),$
there must be a point $\hat{x}$ such that%

\begin{equation}
\hat{x}:\ \ \hat{x}\in\left(  x_{1},x_{2}\right)  \text{ with }v(\hat
{x})=0\text{ where }v\text{ is given by (\ref{v})} \label{vhat}%
\end{equation}
with $v$ increasing at $\hat{x}.$ We further make the following technical
assumption which will be needed for global stability:%
\begin{equation}
\text{The solution to (\ref{vhat})\ is unique with }v^{\prime}\left(  \hat
{x}\right)  >0\text{.} \label{uniq}%
\end{equation}
This holds for a large class of functions and in particular if $f^{\prime
\prime\prime}(x)<0$ for all $x\in\left(  0,a\right)  .$ We identify this point
$\hat{x}$ as the boundary point between mass belonging to spike $\#1$ and mass
belonging to spike $\#2$.

Assuming that the density decays away from the $x_{j},$ we have%
\begin{equation}
M_{1}\sim\int_{-\infty}^{\hat{x}}\rho(x)dx;\ \ \ M_{2}\sim\int_{\hat{x}%
}^{\infty}\rho(x)dx. \label{eqn: def Mj}%
\end{equation}
Integrating (\ref{pde}) from $-\infty$ to $\hat{x}$ we therefore obtain
\begin{equation}
\dfrac{d}{dt}M_{1}=\varepsilon^{2}\rho_{x}(\hat{x}),\ \ \dfrac{d}{dt}%
M_{2}=-\varepsilon^{2}\rho_{x}(\hat{x}) \label{eqn: ODE M}%
\end{equation}
where we used that $v(\hat{x})=0$ and we assumed that $\rho$ decays as
$x\rightarrow\pm\infty.$

To derive an expression for $\rho_{x}(\hat{x})$, next we assume that the
dynamics are sufficiently slow such that $\rho_{t}$ term can be discarded in
(\ref{pde}). This assumption will later be seen to be consisent with the final
result, which shows exponentially slow (metastable)\ evolution of mass. Then
integrating the resulting ODE we obtain%
\begin{equation}
\rho v=\varepsilon^{2}\rho_{x}-\varepsilon^{2}\rho_{x}(\hat{x}),\ \ x\in
\left(  x_{1},x_{2}\right)  , \label{apple}%
\end{equation}
where we use that $v(\hat{x})=0$.

The solution to (\ref{apple}) is given by
\begin{equation}
\rho(x)=\left[  \rho(\hat{x})+\rho_{x}(\hat{x})\,\int_{\hat{x}}^{x}%
\exp(-V(z)/\varepsilon^{2})\,dz\right]  \exp(V(x)/\varepsilon^{2}),
\label{quake}%
\end{equation}
where
\[
V(x)=\int_{\hat{x}}^{x}v(x)dx.
\]

Note that $-V(x)$ has a global maximum at $x=\hat{x}$ (by (\ref{vhat})\ and
(\ref{uniq})), so that we may use Laplace's method to evalute the integral in
(\ref{quake}). We obtain:%
\[
\int_{\hat{x}}^{x}\exp(-V(z)/\varepsilon^{2})\,dz\sim\operatorname*{sign}%
\left(  x-\hat{x}\right)  \varepsilon\sqrt{\frac{\pi}{2v^{\prime}(\hat{x})}%
},\ \ \ \left\vert x-\hat{x}\right\vert \gg O(\varepsilon)
\]
To determine $\rho_{x}\left(  \hat{x}\right)  ,$ let $x\rightarrow
x_{i},\ i=1,2$, and match (\ref{quake})\
with the inner solution as given by (\ref{rhoinner}). Expanding near $x_{j},$
we let $x=x_{j}+\varepsilon y$ and expand%
\[
V(x_{j}+\varepsilon y)=\int_{\hat{x}}^{x_{j}}v(s)ds+v(x_{j})\varepsilon
y+v^{\prime}(x_{j})\varepsilon^{2}\frac{y^{2}}{2}+\ldots\sim\int_{\hat{x}%
}^{x_{j}}v(s)ds-c_{j}\varepsilon^{2}\frac{y^{2}}{2}%
\]
where $c_{j}=-v^{\prime}(x_{j})$ is as given in (\ref{cj}). Therefore the
outer region written in inner variables near $x_{j}$ becomes%
\begin{equation}
\rho(x_{j}+\varepsilon y)\sim\left[  \rho(\hat{x})+\operatorname*{sign}\left(
x_{j}-\hat{x}\right)  \rho_{x}(\hat{x})\,\varepsilon\sqrt{\frac{\pi
}{2v^{\prime}(\hat{x})}}\right]  \exp\left(  \frac{1}{\varepsilon^{2}}%
\int_{\hat{x}}^{x_{j}}v(s)ds\right)  \exp\left(  -c_{j}\frac{y^{2}}{2}\right)
. \label{pomme}%
\end{equation}
On the other hand, the inner region near $x=$ $x_{j}+\varepsilon y$ as derived
in (\ref{rhoinner})\ is
\begin{equation}
\rho\sim\frac{M_{j}}{\varepsilon}\sqrt{\frac{c_{j}}{2\pi}}\exp\left(
-c_{j}\frac{y^{2}}{2}\right)  . \label{inner}%
\end{equation}
Matching (\ref{inner})\ and (\ref{pomme})\ yields\bes\label{oato}%
\begin{align}
\frac{M_{1}}{\varepsilon}\sqrt{\frac{c_{1}}{2\pi}}=\left[  \rho(\hat{x}%
)-\rho_{x}(\hat{x})\,\varepsilon\sqrt{\frac{\pi}{2v^{\prime}(\hat{x})}%
}\right]  \exp\left(  \frac{1}{\varepsilon^{2}}\int_{\hat{x}}^{x_{1}%
}v(s)ds\right),\\
\frac{M_{2}}{\varepsilon}\sqrt{\frac{c_{2}}{2\pi}}=\left[
\rho(\hat{x})+\rho_{x}(\hat{x})\,\varepsilon\sqrt{\frac{\pi}{2v^{\prime}%
(\hat{x})}}\right]  \exp\left(  \frac{1}{\varepsilon^{2}}\int_{\hat{x}}%
^{x_{2}}v(s)ds\right)  .
\end{align}
\ees Solving for $\rho_{x}\left(  \hat{x}\right)  $ from (\ref{oato})\ and
then using (\ref{eqn: ODE M})\ finally yields\bes\label{m1m2}%
\begin{align}
\frac{dM_{1}}{dt}  &  \sim\frac{M_{2}}{2}\sqrt{\frac{c_{2}\,v^{\prime}(\hat
{x})}{\pi^{2}}}\exp\left(  -\frac{1}{\varepsilon^{2}}\int_{\hat{x}}^{x_{2}%
}v(s)ds\right)  -\frac{M_{1}}{2}\sqrt{\frac{c_{1}\,v^{\prime}(\hat{x})}%
{\pi^{2}}}\exp\left(  -\frac{1}{\varepsilon^{2}}\int_{\hat{x}}^{x_{1}%
}v(s)ds\right)  ;\label{M1}\\
\frac{dM_{2}}{dt}  &  =-\frac{dM_{1}}{dt}.
\end{align}
\ees\begin{empty}
\begin{figure}[t]
\centering
\includegraphics[width=0.5\textwidth]{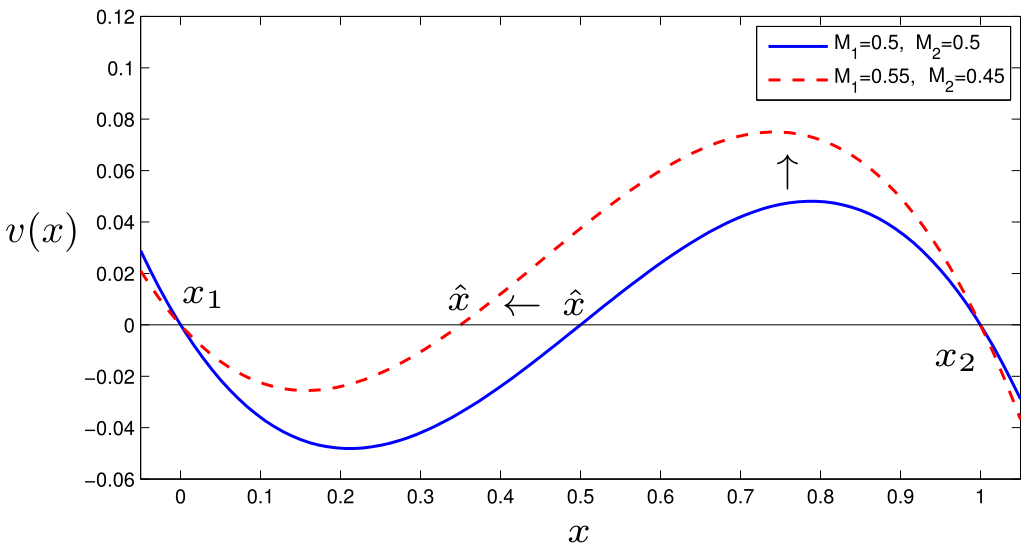}
\caption{The graph of $v(x)=M_1 f(x)+M_2 f(x-1)$ with $f(x)=x-x^3$ and with $M_1,M_2$ as indicated in the
legend. Increasing $M_1$ moves $\hat{x}$ to the left.  }
\label{fig:v}%
\end{figure}
\end{empty}

It is clear that the symmetric configuration $M_{1}=M_{2}$ is an equilibrium
of the ODE\ (\ref{m1m2}), since in this case, $\hat{x}=x_{1}+a/2$ and
$c_{1}=c_{2,}\ \ \int_{\hat{x}}^{x_{1}}v(s)ds=\int_{\hat{x}}^{x_{2}}v(s)ds$.
We now show that it is indeed a global attractor, provided (\ref{uniq}%
)\ holds. From the equation $v(x)\sim M_{1}f(x-x_{1})+M_{2}f(x-x_{2})$ and the
fact that $f(x-x_{1})$ is positive for $x\in\left(  x_{1},x_{2}\right)  ,$ it
follows that $v(x)$ is an increasing function of $M_{1}.$ This in turn shows
that $\hat{x}$ is a decreasing with $M_{1},$ $\int_{\hat{x}}^{x_{1}}v(s)ds$ is
decreasing with $M_{1}$ and $\int_{\hat{x}}^{x_{2}}v(s)ds$ is increasing with
$M_{1}$. Refer to Figure \ref{fig:v}. It then follows from (\ref{m1m2})\ that
$\frac{dM_{1}}{dt}<0$ whenever $M_{1}>M_{2}$ and $\frac{dM_{1}}{dt}>0$ when
$M_{1}<M_{2}.$ This shows that any admissible initial masses satisfying
(\ref{mratio}) evolve towards the equal-mass $M_{1}=M_{2}$ configuration.

We now summarize.

\begin{prop}
\label{prop:ode}Consider the quasi-steady state constructed in Proposition
\ref{prop:ss}. The spike masses $M_{1}(t),\ M_{2}(t)$ evolve on an
exponentially slow time-scale according to (\ref{m1m2}), where $v(x)$ is given
by (\ref{v})\ and $\hat{x}$ satisfies (\ref{vhat}). Moreover, suppose that in
addition to properties (\ref{cond}), $f(x)$ also satisfies: $f^{\prime
\prime\prime}(x)<0$ for $x\in\left(  0,a\right)  .$ Then $M_{1}=M_{2}=M/2$ is
the global attactor of (\ref{m1m2})\ where $M=M_{1}+M_{2}$ is the total mass,
so that $M_{1}(t),M_{2}(t)\rightarrow M/2$ as $t\rightarrow\infty.$
\end{prop}

Note that the expression for $v(s)$ as given by (\ref{v})\ is only accurate up
to $O(\varepsilon^{2})$ terms. As a result, the constants in front of the
exponentials in (\ref{M1})\ are not asymptotically accurate, since the
expansion of the integral in the exponential up to $O(\varepsilon^{2})$ terms
will change the these constants by an $O\left(  1\right)  $ amount.
Nonetheless, the ODE (\ref{m1m2})\ provides a good estimate of the solution at
large time scales as we now show.

To illustrate Proposition \ref{prop:ode}, we take
\[
f(x)=x(1-x^{2}).
\]

Then the right hand side in (\ref{m1m2})\ can be computed explicitly. Without
loss of generality (translation invariance), assume that $x_{1}=0$, hence
$x_{2}=1$ so that $v(x)$ given by (\ref{v})\ becomes
\[
v(x)\sim x(1-x)\left(  M_{1}-2M_{2}+M\,x\right)  .
\]
The unique $\hat{x}\in(0,1)$ for which $v(\hat{x})=0$ is therefore explicitly
given by
\[
\hat{x}=\dfrac{2M_{2}-M_{1}}{M},
\]
and we compute
\begin{equation}
v^{\prime}(\hat{x})=\,\dfrac{(2M_{2}-M_{1})(2M_{1}-M_{2})}{M}; \label{1003}%
\end{equation}
and%
\begin{equation}
\int_{\hat{x}}^{x_{1}}v\left(  x\right)  dx=\dfrac{M_{1}(2M_{2}-M_{1})^{3}%
}{4M^{3}};\ \ \ \ \int_{\hat{x}}^{x_{2}}v\left(  x\right)  dx=\dfrac
{M_{2}(2M_{1}-M_{2})^{3}}{4M^{3}}. \label{1004}%
\end{equation}
Substituting (\ref{1003})\ and (\ref{1004})\ into (\ref{m1m2}) yields the
explicit ODE\ (\ref{odecubic})\ for the mass exchange dynamics.

As mentioned above, the constants in front of exponentials are not
asymptotically accurate, since we only computed $O(1)$ contributions of the
arguments $\int_{\hat{x}}^{x_{1}}v\left(  x\right)  dx,\ \int_{\hat{x}}%
^{x_{2}}v\left(  x\right)  dx$ inside the exponents, disregarding the
$O(\varepsilon^{2})$ corrections which contribute an $O(1)$ amount to the
constants in front of the exponential. However the ODE is still asymptotically
valid for large time in the following sense. Rewrite the system
(\ref{odecubic}) in terms of the difference
\[
d=M_{2}-M_{1}%
\]
from the equilibrium state. If we assume that $d$ is positive and $d\gg
O(\varepsilon^{2})$ then the term $F(M_{2},M_{1})$ in (\ref{odeM})$\ $is
negligible when compared with $F(M_{1},M_{2})$ and we obtain%
\begin{equation}
d^{\prime}(t)\sim-C(d)\exp\left\{  -\frac{1}{\varepsilon^{2}}\frac{\left(
M+d\right)  \left(  M-3d\right)  ^{3}}{64M^{3}}\right\}  \label{dp}%
\end{equation}
where $C(d)$ is some (undetermined) $O(1)$ pre-multiplier. Taking the logarithm on
both sides of (\ref{dp})\ we obtain
\begin{equation}
\varepsilon^{2}\log\left(  \left\vert d^{\prime}(t)\right\vert \right)
=\frac{-\left(  M+d\right)  \left(  M-3d\right)  ^{3}}{64M^{3}}+O(\varepsilon
^{2}). \label{logdp}%
\end{equation}
Figure \ref{fig:ddp}(left) shows direct comparison between the full numerical
simulations of the original PDE\ (\ref{pde})\ and the asymptotics of the right-hand side given by (\ref{logdp}). The dashed asymptotic line there corresponds
to the right-hand side of (\ref{logdp}). Note how the difference between the
dashed line (asymptotics)\ and solid lines\ (numerics for several
$\varepsilon$)\ is proportional to $\varepsilon^{2},$ which shows that the
error in the exponent in (\ref{M1})\ is indeed of $O(\varepsilon^{2}).$

The ODE\ (\ref{dp})\ is separable, and applying Laplace's method to its
solution (keeping in mind that $\left(  M+d\right)  \left(  M-3d\right)  ^{3}$
is a decreasing function of $d$) we obtain the asymptotic solution
\begin{equation}
\left(  M+d\right)  \left(  M-3d\right)  ^{3}\sim64M^{3}\varepsilon^{2}\left(
\log\left(  \frac{t}{\varepsilon^{2}64M^{3}}\right)  +O(1)\right)
\label{logt}%
\end{equation}
valid for large $t.$ Note that to leading order, $d(t)$ is independent of the
pre-exponential multiplier $C(d)$ in\ (\ref{dp}) (it appears in the $O(1)$
term in the right hand side), and is also indepedent of the initial conditions
as long as $d\gg O(\varepsilon^{2})$. In other words (\ref{logt}) is an
accurate prediction for $d$ as a function of $\log t$, regardless of the value
of $C(d).$

\begin{empty}
\begin{figure}[t]
\includegraphics[width=0.49\textwidth]{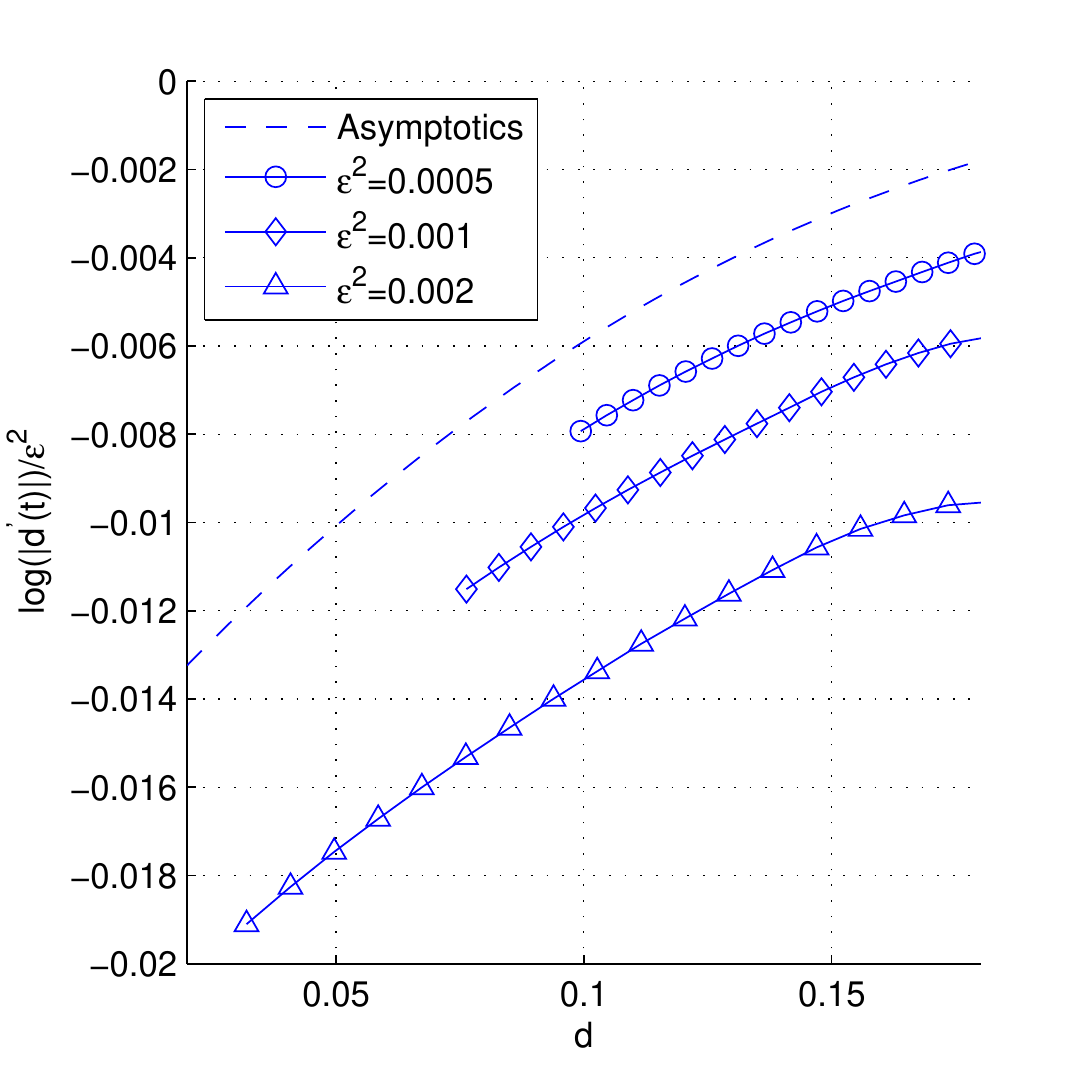}
\includegraphics[width=0.49\textwidth]{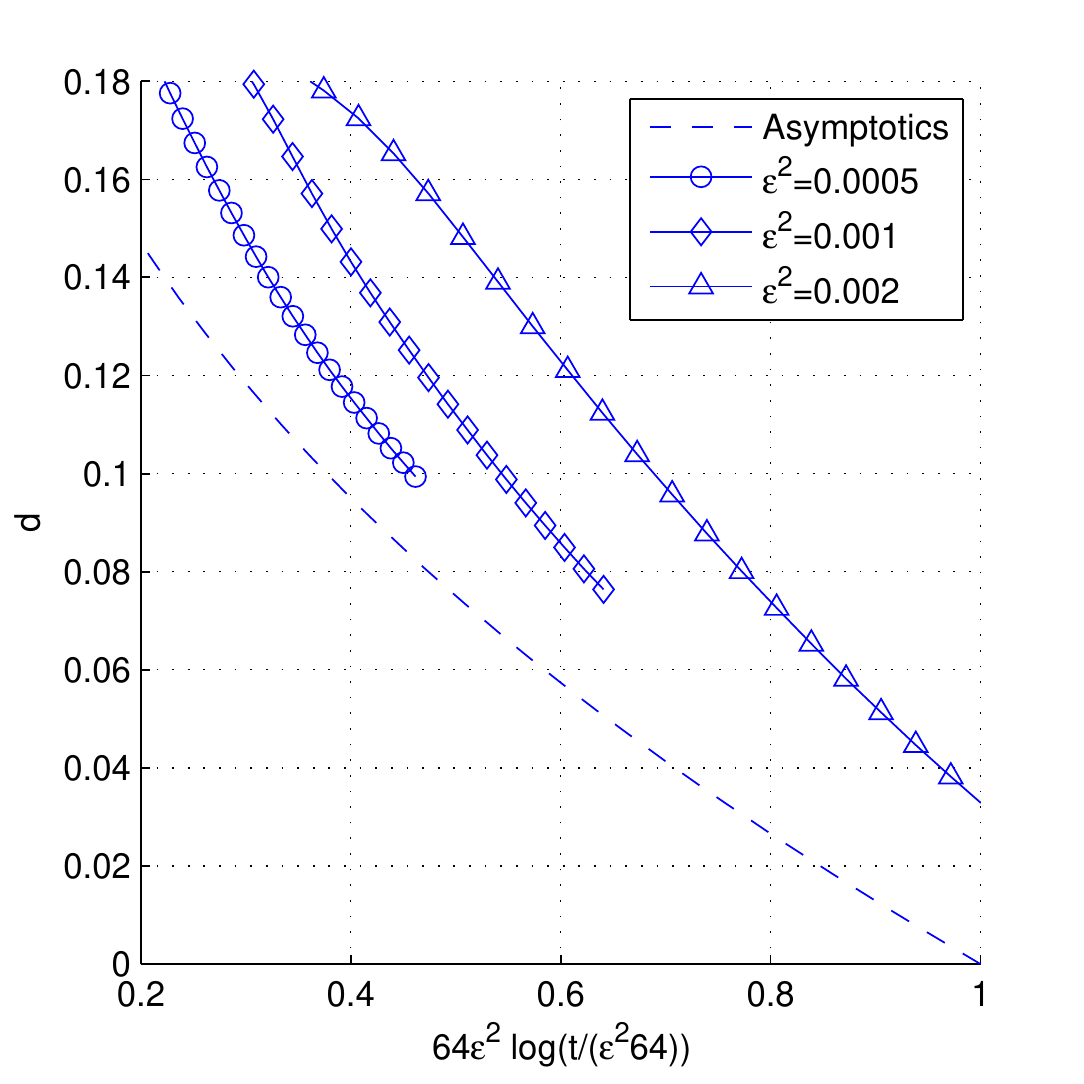}
\caption{LEFT:
Validation of the right hand side of the ODE (\ref{logdp}). 
The plot shows $d(t)$ versus ${\varepsilon^2}\log(|d'(t)|)$, where $d(t)$
and $d'(t)$ are computed from numerical simulations of the full PDE
for several values of $\varepsilon$ as shown.
The dashed line is the leading-order
asymptotic prediction (\ref{logdp}) ${\varepsilon^2}\log(|d'(t)|)\approx
-\frac{(M+d)(M-3d)^3}{64M^3}$.
RIGHT: Same simulation/parameters as on the left, except that the
mass difference $d(t)$ is plotted versus logarithmic timescale. The theoretical
prediction plotted in dashed line is given by (\ref{logt}).
}
\label{fig:ddp}%
\end{figure}
\end{empty}

In Figure \ref{fig:ddp}(right) we plot $d$ as a function of $64M^{3}%
\varepsilon^{2}\left(  \log\left(  \frac{t}{\varepsilon^{2}64M^{3}}\right)
\right)  $ for several values of $\varepsilon,$ and then compare with the
asymptotic value of $d$ as given by\ (\ref{logt}). Again, the error is
observed to be proportional to $\varepsilon^{2}.$

For the full numerical solution to the original PDE
(\ref{pde}), we took the computational domain to be
of size $x\in\left[  0,3\right]$. Because of the exponential decay outside the spikes, doubling the domain size did not change the results. The initial
conditions were taken to be (\ref{rhoinner}), with $M_{1}=0.35,\ M_{2}=0.65$ and the centers of the two Gaussian peaks were a distance $1$ apart. We also
waited $t=10$ time units to let the transients die out before starting the
comparison; at $t=10$ the system already converged to the quasi-steady state. The masses $M_{1}$
and $M_{2}$ in (\ref{v})\ were computed numerically using (\ref{eqn: def Mj})
for each time-step. We used finite
differences with semi-implicit time stepping:\ $v(x)$ is computed explicitly
at each time step, while the update for $\rho(x,t+\Delta t)$ is done
implicitly. We verified the accuracy by using several stepsizes.

\section{Discussion}

\label{sec:discuss}

In this paper we focused on the steady states of\ the aggregation equation
with noise (\ref{pde}) that consist of two nearly-Dirac concentrations. An
important implication of the work presented here, is that not all steady
states of the zero-diffusion equation can be recovered as the limit as
$\varepsilon\rightarrow0$ of a sequence of steady states of the
aggregation-diffusion equation. This concerns in particular the two-Dirac
steady states of unequal mass. Nevertheless, on short timescales some
reminiscents of these unequal-mass steady states are still present. To be more
precise, there is an intermediate timescale in which a state consisting of two
Diracs of unequal mass persists as a metastable state.

We have shown that the process of equilibration can be described
asymptotically by an ODE for the evolution of the mass associated to each of
the spikes. This ODE is valid on an exponentially long timescale, after the
initial two-spike quasi-steady state profile is formed on an $O(1)$ timescale.
An open question is to derive a more precise ODE for the dynamics, which would
require the next-order expansion of the quasi-steady state; see figure
\ref{fig:ddp} and associated discussion in \S \ref{sec: ODE derivation}.

An interesting, but nontrivial complication arises when one wants to derive an
ODE for the mass evolution for three or more spikes. In the case of two
spikes, these are centred a distance $a$ apart (with $a$ such that $f(a)=0$).
For more than two spikes a delicate balance needs to be satisfied: there is an
algebraic system of equations involving the masses and mutual distances, such
that the velocity at each centre is zero. This algebraic relation, together
with an ODE for the evolution of each mass, accounts for the simultaneous
evolution of the centres of the spikes and the masses towards equilibrium.

As in the case of two spikes, the steady state consisting of three
delta-concentrations in the absence of diffusion is degenerate: there is an
arbitrariness in how the three masses are distributed among the three holes,
and there are two degrees of freedom (three masses subject to constraint
$M_{1}+M_{2}+M_{3}=1$). However when the diffusion is turned on, this
two-parameter family of steady states should ``collapse'' into a unique steady
state. It would be very interesting to characterize precisely which mass
fractions are ``selected'' by the diffusion. Unlike the two-spike solution
where diffusion ``chooses'' the equal-mass configuration, in the case of the
three-spike configuration, the diffusion should in general select unequal mass fractions.

It would be interesting to extend these results to two and higher dimensions.
In two dimensions, at least three delta-concentrations are required for
stability \cite{kolokolnikov2013singular} (in the case three
delta-concentrations, their locations form an equilateral triangle). The
construction of the inner solution near the spike is analogous to the
derivation in \S \ref{sec:ss}. On the other hand, the outer region cannot be
easily solved, as it requires solving a fully two-dimensional PDE, and
performing matching is a nontrivial problem. Nonetheless there is a hope that
WKB-type techniques can be used to approximate the solution to the outer
region for small diffusion. It would be interesting if similar
\textquotedblleft equilibration\textquotedblright\ results can be obtained in
two dimensions.

\section{Acknowledgements}

We thank the anonymous referees for their careful reading of the first draft
of the paper and for pointing out that the constant $C(d)$ in the ODE
(\ref{dp})\ was not fully resolved. J.E.~is supported by an AARMS Postdoctoral
Fellowship. T.K.~is supported by NSERC discovery and NSERC
accelerator grants.

\bibliographystyle{elsarticle-num}
\bibliography{2masses}

\end{document}